\documentclass{elsarticle}      

\usepackage{graphicx}

\begin{document}

\markboth{P. Kalinay, L. \v{S}amaj, I. Trav\v{e}nec}
{Survival probability (heat content) and the lowest eigenvalue 
of Dirichlet Laplacian}


\title{Survival probability (heat content) and the lowest eigenvalue 
of Dirichlet Laplacian}

\author{Pavol Kalinay\footnote{Pavol.Kalinay@savba.sk},
Ladislav \v{S}amaj, Igor Trav\v{e}nec}
\address{Institute of Physics, Slovak Academy of Sciences, \\
D\'ubravsk\'a cesta 9, 845 11 Bratislava, Slovakia}

\maketitle


\begin{abstract}
We study the survival probability of a particle diffusing in a
two-dimensional domain, bounded by a smooth absorbing boundary.
The short-time expansion of this quantity depends on the geometric 
characteristics of the boundary, whilst its long-time asymptotics is 
governed by the lowest eigenvalue of the Dirichlet Laplacian defined on 
the domain. 
We present a simple algorithm for calculation of the short-time
expansion for an arbitrary "star-shaped" domain.
The coefficients are expressed in terms of powers of boundary curvature,
integrated around the circumference of the domain. 
Based on this expansion, we look for a Pad\'e interpolation between 
the short-time and the long-time behavior of the survival probability, 
i.e. between geometric characteristics of the boundary and
the lowest eigenvalue of the Dirichlet Laplacian.
\end{abstract}


\section{Introduction}
Many problems in physics require to study how the overall properties of
a system are influenced by its boundary. 
An archetypal task is looking for dependence of the lowest eigenvalue 
of the (minus) Laplacian, defined in a finite domain, 
on the shape of this domain. 
This problem is related to topics like the calculation of 
the ground state energy of a free quantum particle in a finite domain, 
the lowest cutoff frequency of the dominant mode in a waveguide, 
or the lowest tone of a drum\cite{Kac66,Prot87}. 
We shall formulate this mathematical problem in the context of the
diffusion (probabilistic) theory, which is exploited often in 
the chemical physics: The geometry of a molecular boundary influences 
the reaction rates\cite{Sano79}, or the dynamics of particles 
diffusing in confined systems\cite{Mon06,Kalinay07,Kalinay08}.
Our attention is focused on the survival probability of a particle
diffusing in a finite domain, which is defined later.
We show that this quantity interrelates the boundary characteristics and
the lowest eigenvalue of the Dirichlet Laplacian in an interesting way. 

To introduce the notation, we consider a finite domain $\Omega$ 
of points ${\bf r}$ with a smooth boundary $\partial\Omega$.
The spectrum of the (minus) Laplacian, say with the Dirichlet boundary
condition (BC), is given by
\begin{equation} \label{e1}
\begin{array}{rcll}
-\Delta \phi({\bf r}) & = & \lambda \phi({\bf r}) 
&\qquad {\bf r}\in \Omega , \cr
\phi({\bf r}) & = & 0 &\qquad {\bf r}\in \partial\Omega .
\end{array}
\end{equation}
The eigenvalues $0<\lambda_1<\lambda_2\le\lambda_3\cdots\le\lambda_j\le\cdots$
form a discrete set\cite{Courant53}.
The corresponding eigenfunctions $\phi_1,\phi_2,\phi_3,\ldots,\phi_j,\ldots$
form an orthonormal basis of real functions,
\begin{equation} \label{e2}
\int_{\Omega} {\rm d}{\bf r} \phi_j({\bf r}) \phi_k({\bf r}) = \delta_{jk} ; 
\end{equation}
the eigenfunctions satisfy the completeness relation
\begin{equation} \label{e3}
\sum_j \phi_j({\bf r}) \phi_j({\bf r}') = \delta({\bf r}-{\bf r}') .
\end{equation} 

In the diffusion theory, the conditional probability 
$\rho({\bf r},t\vert {\bf r}_0,0)$ of finding a particle at a point 
${\bf r} \in \Omega$ at time $t>0$, if it started from ${\bf r}_0\in \Omega$ 
at $t_0=0$, is governed by the diffusion (heat) equation
\begin{equation} \label{e4}
\frac{\partial \rho({\bf r},t\vert {\bf r}_0,0)}{\partial t} =
\Delta \rho({\bf r},t\vert {\bf r}_0,0) ,
\end{equation}
where the diffusion constant is set to 1.
This equation has to be supplemented by the initial condition
$\rho({\bf r},t=0\vert {\bf r}_0,0) = \delta({\bf r}-{\bf r}_0)$
and by the Dirichlet BC $\rho({\bf r},t\vert {\bf r}_0,0) = 0$ for 
${\bf r}\in\partial\Omega$, which reflects absorption (disappearance) 
of the particle hitting the boundary.
The conditional probability can be expressed in terms of the eigenvalues
and eigenfunctions of the Dirichlet Laplacian as follows
\begin{equation} \label{e5}
\rho({\bf r},t\vert {\bf r}_0,0) = \sum_j \phi_j({\bf r}_0) \phi_j({\bf r})
{\rm e}^{-\lambda_j t} .
\end{equation}

The ``local'' survival probability is defined as
\begin{equation} \label{e6}
S(t,{\bf r}_0) = \int_{\Omega} {\rm d}{\bf r} 
\rho({\bf r},t\vert {\bf r}_0,0) .
\end{equation}
It represents the probability that a particle, localized at a point
${\bf r}_0\in\Omega$ at the initial time $t_0=0$, remains still diffusing
in the domain $\Omega$ at time $t>0$ unabsorbed by the boundary.
Since the particle inserted at the boundary at $t_0=0$ is immediately 
absorbed, $S(t,{\bf r}_0)$ satisfies also the Dirichlet BC
\begin{equation} \label{e7}
S(t,{\bf r}_0) = 0 \qquad \mbox{for ${\bf r}_0\in \partial\Omega$.}
\end{equation} 
The quantity of our interest is the ``global'' survival probability $S(t)$, 
defined as the average of the local one over the whole domain,
\begin{equation} \label{e8}
S(t) = \frac{1}{\vert\Omega\vert} \int_{\Omega} {\rm d}{\bf r}_0
S(t,{\bf r}_0) = \sum_{j=1}^{\infty} \gamma_j^2 e^{-\lambda_jt} , \qquad
\gamma_j = \frac{1}{\sqrt{\vert\Omega\vert}} 
\int_{\Omega} {\rm d} {\bf r} \phi_j({\bf r}) .
\end{equation}
It represents the probability of finding the particle in $\Omega$ at time 
$t>0$, if it was distributed uniformly with the density $1/\vert\Omega\vert$
over the whole domain at $t_0=0$. 

The same analysis can be applied to the heat contained in the domain
$\Omega$, heated at some nonzero temperature, whose boundary 
$\partial\Omega$ is cooled to the zero temperature for $t>0$. 
$S(t)$ is then proportional to {\it the heat content} in the domain and 
this is the name under which this quantity is known in the mathematical 
literature\cite{Birkhoff54,Berg93,Berg94,DesJardins98,Savo98a,Savo98b}.

First let us review some basic properties of $S(t)$.
At $t=0$, it starts from $S(0)=1$ as the particle just inserted
somewhere in the domain is surely not yet absorbed by the boundary;
the consequent equality $\sum_j \gamma_j^2 =1$ follows trivially
from the completeness relation (\ref{e3}).
The monotonous decay of $S(t)$ at small $t$ is determined by the geometry of 
the boundary, as only the particle inserted close to the boundary have
a chance to be absorbed very soon.
Like for any spectral function, the small-$t$ expansion of $S(t)$ can
be obtain explicitly.
It is an expansion in $\sqrt{t}$ of the form  
\begin{equation}\label{e9}
S(t) = 1 + \sum_{j=1}^{\infty} \sigma_j t^{j/2} ,
\end{equation}
for any dimensionality of the domain\cite{Savo98a,Savo98b}.
The coefficients $\{ \sigma_j\}$ were shown to satisfy a recurrence scheme, 
containing certain differential operators acting on the shortest distance 
$l({\bf r})$ of points ${\bf r}$ from the boundary
$\partial\Omega$\cite{Savo98a,Savo98b}. 
The forms of the coefficients become soon very complicated,
only first few of them are available explicitly.
The first aim of this paper is to approximate reasonably the coefficients
$\{ \sigma_j\}$ in higher orders by local curvatures of the domain boundary. 

It is clear that the whole eigenspectrum contributes to the coefficients 
$\{ \sigma_j \}$. 
On the other hand, according to (\ref{e8}), the long-time decay of 
the survival probability to 0 is governed by the lowest eigenvalue of
the Dirichlet Laplacian $\lambda_1$, 
\begin{equation} \label{e10}
S(t)\sim\gamma_1^2 \exp(-\lambda_1 t) \qquad \mbox{for $t\to\infty$}.
\end{equation}
Thus the survival probability interconnects the geometric characteristics of 
the domain boundary, contained in the small-$t$ expansion of $S(t)$, with 
the overall characteristics of the system, represented by the lowest eigenvalue
$\lambda_1$ in the long-time limit. 
Our second aim is to apply some appropriate interpolation scheme 
between the small-$t$ and large-$t$ expansions.
In this way we can investigate how the lowest eigenvalue of the Dirichlet
Laplacian in a domain depends on the geometry of its boundary. 

We are inspired by the method of t-expansions for studying
the ground state energy of many-body systems\cite{Horn84}. 
The method is based on an energy functional $E(t)$, depending on time $t$.
This functional has an analytic small-$t$ expansion in terms of 
connected moments $\{ I_j \}$, $E(t)=\sum_{j=0}^{\infty} (-t)^j I_{j+1}/j!$. 
The ground state energy $\lambda_1$ follows from the long time behavior 
$\lambda_1=\lim_{t\to\infty} E(t)$.
There exist several interpolation schemes between the small-$t$ and large-$t$
expansions\cite{Cioslowski87a,Cioslowski87b,Stubbins88,Samaj97} 
which enable one to calculate $\lambda_1$ in terms of connected moments 
$\{ I_j \}$.
The interpolation schemes are based on various combinations of exponential
decays of $E(t)$ to the asymptotic value $\lambda_1$; the fact that
the small-t expansion is analytic plays a fundamental role.
The method of t-expansions was successfully used for calculation of the
ground state energy of 
molecules\cite{Fessatidis02,Fessatidis05,Fessatidis07,Fessatidis08}, or for 
the study of phase transitions in infinite lattice systems\cite{Samaj97}.
Since the small-$t$ expansion of the survival probability (\ref{e9}) is
in powers of $\sqrt{t}$, the interpolation schemes working in $t$-expansions 
are not straightforwardly applicable to our problem. 

As was indicated, this paper is devoted to and brings new results in two 
topics:
\begin{itemize}
\item
A reasonable approximation of the coefficients $\{ \sigma_j\}$ of 
the small-$t$ expansion (\ref{e9}) in high orders.
\item
The search for an interpolation formula between the small-$t$ (\ref{e9}) 
and large-$t$ (\ref{e10}) expansions of $S(t)$.
\end{itemize}
It turns out to be advantageous to perform the analysis in the Laplace space.
In Sect. 2, we introduce the Laplace transform for the local and global
survival probabilities and derive a differential equation for the former.
In Sect. 3, we use the exact solution for a disk domain to approximate
reasonably the coefficients $\{ \sigma_j\}$ of the small-$t$ expansion 
(\ref{e9}) up to high orders which is desirable for the proposed
interpolation method. 
The first four coefficients are identical to the exact 
formulas\cite{Savo98a,Savo98b}, the next approximative ones involve correctly
the powers of the boundary curvature integrated along $\partial\Omega$. 
Sect. 4 is devoted to an interpolation method between the small-$t$
and large-$t$ expansions of $S(t)$, based on a variant of the Pad\'e 
approximant of the Laplace transform of $S(t)$.
Our algorithm is tested on the circular and elliptical domains.

\section{Laplace transform of $S(t)$}
Our method is formulated in the Laplace space.
We introduce the Laplace transforms with respect to $t$, with $s^2$ as 
the Laplace parameter, for both the local survival probability
\begin{equation} \label{e13}
\tau(s,{\bf r}_0) = \int_0^{\infty} {\rm d}t {\rm e}^{-s^2 t} S(t,{\bf r}_0) ,
\qquad \tau(s,{\bf r}_0) = 0 \quad \mbox{for ${\bf r}_0\in\partial\Omega$}
\end{equation}
and the global survival probability
\begin{equation} \label{e14}
\tau(s) = \int_0^{\infty} {\rm d}t {\rm e}^{-s^2 t} S(t)
= \sum_{j=1}^{\infty} \frac{\gamma_j^2}{s^2+\lambda_j} . 
\end{equation}
The Laplace transforms are related by
\begin{equation} \label{e15}
\tau(s) = \frac{1}{\vert\Omega\vert} \int_{\Omega} {\rm d}{\bf r}_0
\tau(s,{\bf r}_0) .
\end{equation}
In the Laplace picture, the small-$t$ behavior of $S(t)$ (\ref{e9})
is represented as a large-$s$ expansion of $\tau(s)$ in $1/s$:
\begin{equation} \label{e16}
\tau(s) = \frac{1}{s^2} + \sum_{j=1}^{\infty}
\frac{\Gamma(j/2+1)}{s^{j+2}}\sigma_j .
\end{equation}
The large-$t$ behavior of $S(t)$ is represented as a regular small-$s$
expansion of $\tau(s)$ in $s^2$:
\begin{equation} \label{e17}
\tau(s) = \sum_{k=0}^{\infty} (-1)^k s^{2k} \sum_{j=1}^{\infty}
\frac{\gamma_j^2}{\lambda_j^k} ,
\end{equation} 
with a dominant contribution of the lowest eigenvalues.

The Laplace transform of the local survival probability satisfies 
a differential equation which is derived by the following 
procedure\cite{Zwanzig01}.
If we apply the conjugated Laplacian $\Delta^+$ (acting upon the
coordinates of ${\bf r}_0$) to both sides of (\ref{e13}) and use 
the conjugated of the diffusion equation (\ref{e4}), we obtain  
\begin{equation} \label{e18}
\Delta^+ \tau(s,{\bf r}_0) = \int_0^{\infty} {\rm d}t {\rm e}^{-s^2 t}
\int_{\Omega} {\rm d}{\bf r} \Delta^+ \rho({\bf r},t\vert {\bf r}_0,0)
= \int_{\Omega} {\rm d}{\bf r} \int_0^{\infty} {\rm d}t
{\rm e}^{-s^2 t} \partial_t\rho({\bf r},t\vert {\bf r}_0,0) .
\end{equation}
Integration by parts in $t$ then implies
\begin{equation} \label{e19}
\Delta \tau(s,{\bf r}) - s^2 \tau(s,{\bf r}) = - 1 ,
\end{equation}
where we abandon the subscript 0 and replace $\Delta^+$ by $\Delta$.
This differential equation is supplemented by the Dirichlet BC
\begin{equation} \label{e20}
\tau(s,{\bf r}) = 0 \qquad \mbox{for ${\bf r}\in \partial\Omega$.}
\end{equation}

\begin{figure}
\begin{center}
\includegraphics[scale=0.45]{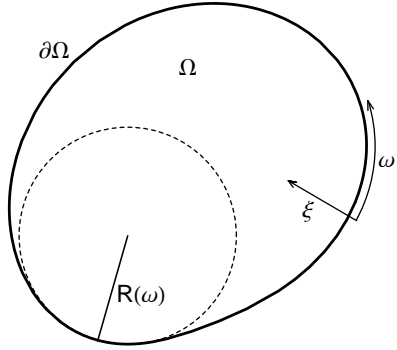}
\end{center}
\caption{Fermi coordinates $(\omega,\xi)$ for general domain $\Omega$ 
with boundary $\partial\Omega$.}
\end{figure}

In the remaining part of the paper, we shall restrict ourselves to 
two-dimensional star-shaped domains with only {\em one} smooth boundary.
To express $\tau(s)$ for such domains, it is convenient to average 
the differential equation (\ref{e19}) over $\Omega$ and then to apply 
the Gauss-Ostrogradsky theorem:
\begin{equation} \label{e21}
s^2 \tau(s) = 1 + \frac{1}{\vert\Omega\vert} \int_{\Omega} {\rm d}{\bf r} 
\Delta\tau(s,{\bf r}) = 1 - \frac{1}{\vert\Omega\vert}
\int_{\partial\Omega} {\rm d}\omega 
\frac{\partial\tau(s;\omega,\xi)}{\partial \xi}\Big\vert_{\xi=0} .
\end{equation}
Here, we introduced the Fermi coordinates\cite{Chavel84} (see Fig. 1);
$\omega$ runs around the domain boundary and $\xi$ is orthogonal to $\omega$ 
such that $\xi=0$ at the boundary $\partial\Omega$ and $\xi>0$ inside
the domain $\Omega$.

\section{Approximation of $S(t)$ by local curvature}

\subsection{Disk domain}
The simplest star-shaped domain is a disk of radius $R$, with
the boundary defined by $r(\varphi)=R$.
The disk geometry is the crucial one in our approach. 
The eigenfunctions of the Dirichlet Laplacian for the disk are expressed 
in polar coordinates as follows
\begin{equation} \label{e22}
\phi_{n,l}(r,\varphi) = c_{n,l} J_l(r\sqrt{\lambda_{n,l}}) \cos(l\varphi) ,
\end{equation}
where $n=1,2,\ldots$ and $l=0,1,2,\ldots$ are the radial and angular quantum 
numbers, respectively, $c_{n,l}$ are the normalization constants and 
the eigenvalues $\lambda_{n,l}=(z_{l,n}/R)^2$ are determined by the zeros 
$z_{l,n}$ of the Bessel functions of the first kind $J_l(z_{l,n})=0$.
The lowest eigenvalue $\lambda_1\equiv \lambda_{1,0} = z_{0,1}^2/R^2 
= 5.783196\ldots/R^2$.
Only the radial modes
$\phi_{n,0}(r,\varphi)={J_0(z_{0,n}r/R)/\sqrt{\pi} R J_1(z_{0,n})}$
give nonzero integrals 
$\int_{\Omega} r {\rm d}r {\rm d}\varphi \phi_{n,0}(r,\varphi)=
2\sqrt{\pi}R/ z_{0,n}$
in equation (\ref{e7}) and therefore contribute to $S(t)$, i.e.
\begin{equation}\label{e23}
S(t) = 4 \sum_{n=1}^{\infty} \frac{1}{z_{0,n}^2}
\exp\left( - \frac{z_{0,n}^2}{R^2} t \right) .
\end{equation}
This means that the survival probability does not carry information
about the whole spectrum of the Dirichlet Laplacian, only the eigenstates 
with the angular quantum number $l=0$ are involved.
This is a specific property of the disk domain. 

The representation (\ref{e23}), written in terms of the zeros of 
the Bessel function $J_0$, is not suitable for our purposes.
Instead, we use the differential equation (\ref{e19}) supplemented by 
the Dirichlet BC (\ref{e20}).
In polar coordinates, the general solution of (\ref{e19}), 
regular at $r\to 0$, reads
\begin{equation} \label{e24}
\tau(s; r,\varphi) = \frac{1}{s^2} + \sum_{l=-\infty}^{\infty}
\tau_l I_l(s r) e^{i l\varphi} ,  
\end{equation}
where $\{ I_l\}$ are the modified Bessel functions of the first kind.
Since $\tau(s; r,\varphi)$ is real, the complex coefficients 
$\{ \tau_l\}$ must satisfy $\tau_{-l}=\tau_l^*$.
The boundary condition $\tau(s;r=R,\varphi)=0$ implies
\begin{equation} \label{e25}
\tau_0 = - \frac{1}{s^2} \frac{1}{I_0(sR)} , \qquad
\mbox{$\tau_l = 0$ for $l\ne 0$.}
\end{equation}
Consequently,
\begin{equation} \label{e26}
\tau(s; r,\varphi) \equiv \tau(s,r)
= \frac{1}{s^2} \left[ 1 - \frac{I_0(sr)}{I_0(sR)} \right] .
\end{equation}
The averaging of this relation over the disk area leads to the explicit 
form of the Laplace transform of the survival probability\cite{Berg94}
\begin{equation} \label{e27}
\tau(s) = \frac{1}{s^2} \left[ 1- \frac{2 I_1(sR)}{s R I_0(sR)} \right] .
\end{equation}

To obtain the small-$s$ and the large-$s$ expansions of $\tau(s)$, 
we use the corresponding expansions of the Bessel functions 
$I_0(sR)$ and $I_1(sR)$\cite{GR}.
The small-$s$ expansion becomes
\begin{equation} \label{e28}
\tau(s) = \frac{R^2}{8} - \frac{s^2 R^4}{48} + \frac{11 s^4 R^6}{3072}
+ O(s^6) .
\end{equation}
In the limit $s\to\infty$, using the asymptotic expansions of the
Bessel functions $I_0$ and $I_1$, we obtain
\begin{equation} \label{e29}  
s^2 \tau(s) = 1 - \frac{2}{sR} \left[ 1 - \frac{1}{2sR} - 
\frac{1}{8(sR)^2} \right. \left. 
- \frac{1}{8(sR)^3} - \frac{25}{128(sR)^4} - \cdots \right] .
\end{equation}
Comparing the series representations (\ref{e9}) and (\ref{e16}),  
we get the corresponding small-$t$ expansion of the survival probability:
\begin{equation} \label{e30}
S(t) = 1 - \frac{4\sqrt{t}}{\sqrt{\pi}R} + \frac{t}{R^2}
+ \frac{t^{3/2}}{3\sqrt{\pi}R^3} 
+ \frac{t^2}{8 R^4} + \frac{5 t^{5/2}}{24\sqrt{\pi} R^5} + \cdots .
\end{equation}

\subsection{Convex domain}
Being motivated by the exact solution for the disk, we now propose an 
approximate form of the large-$s$ expansion of $\tau(s)$ 
for an arbitrary star-shaped domain. 
In the limit $s\to\infty$, $\tau(s,{\bf r})\simeq 1/s^2$ inside 
the domain and it decreases to 0, given by the Dirichlet BC, only in 
a thin layer near the boundary.
This can be seen from the result (\ref{e26}) if we use the asymptotic 
expansion $I_0(sr)/I_0(sR) \sim e^{-s(R-r)}$ for $s\to\infty$.
Our strategy is to replace the adjacent part of the boundary 
adequately by some simpler curve, for which the function 
$\tau(s; {\bf r})$ is known, and use this function instead of 
the exact one in the formula (\ref{e21}).

The simplest way is to approximate the convex parts of the boundary by circles 
of radius $R(\omega)$, inverse to the curvature of the boundary at 
a given point $\omega$ (see Fig. 1).
Then, instead of the true $\tau(s;\omega,\xi)$, we adopt the disk solution 
(\ref{e26}) with $r$ replaced by $r=R(\omega)-\xi$, so that
\begin{equation}\label{e31}
\frac{\partial\tau(s; \omega,\xi)}{\partial\xi}\Big\vert_{\xi=0}
\simeq \frac{I_1(sR(\omega))}{s I_0(sR(\omega))} .
\end{equation}
Applying the asymptotic $s\to\infty$ expansions of the Bessel functions in 
this formula, the relation (\ref{e21}) can be rewritten directly as 
an expansion in $1/s$:
\begin{eqnarray}
s^2\tilde{\tau}(s) & = & 1 - \frac{1}{s\vert\Omega\vert} \int_{\partial\Omega} 
{\rm d}\omega \left[ 1- \frac{1}{2sR(\omega)} \right. 
- \frac{1}{8 s^2 R^2(\omega)} - \frac{1}{8 s^3 R^3(\omega)} \nonumber \\ & & 
\left. 
- \frac{25}{128 s^4 R^4(\omega)} - \frac{13}{32 s^5 R^5(\omega)} 
\cdots \right] . \label{e32}
\end{eqnarray}
Tilde on the top of $\tilde{\tau}(s)$ means the approximation 
by the local curvature. 
According to the representation (\ref{e16}), the corresponding coefficients 
$\{ \tilde{\sigma}_j \}$ are given by
\begin{eqnarray}\label{33}
\tilde{\sigma}_1 & = &  - \frac{2}{\sqrt{\pi}} \frac{1}{\vert\Omega\vert} 
\int_{\partial\Omega} {\rm d}\omega , \label{e33a} \\
\tilde{\sigma}_2 & = &  \frac{1}{2} \frac{1}{\vert\Omega\vert}
\int_{\partial\Omega} {\rm d} \omega k(\omega) , \label{e33b} \\
\tilde{\sigma}_3 & = & \frac{1}{6\sqrt{\pi}} \frac{1}{\vert\Omega\vert}
\int_{\partial\Omega} {\rm d}\omega k^2(\omega) , \label{e33c} \\
\tilde{\sigma}_4 & = &  \frac{1}{16} \frac{1}{\vert\Omega\vert}
\int_{\partial\Omega} {\rm d} \omega k^3(\omega) , \label{e33d} \\
\tilde{\sigma}_5 & = &  \frac{5}{48\sqrt{\pi}} \frac{1}{\vert\Omega\vert}
\int_{\partial\Omega} {\rm d} \omega k^4(\omega) , \label{e33e} \\
\tilde{\sigma}_6 & = &  \frac{13}{192} \frac{1}{\vert\Omega\vert}
\int_{\partial\Omega} {\rm d} \omega k^5(\omega) , \label{e33f}
\end{eqnarray}
etc., where $k(\omega)\equiv 1/R(\omega)$ is the curvature at the point 
$\omega\in \partial\Omega$.
It follows from the construction that the coefficients can be expressed
explicitly to very high orders.
If the boundary is defined in polar coordinates,
$r=r(\varphi)$ with $\varphi\in [0,2\pi]$,
the length element of the boundary $\partial\Omega$ is given by
\begin{equation} \label{e34}
d\omega = \sqrt{r^2+r'^{2}} d\varphi
\end{equation}
and the curvature by
\begin{equation} \label{e35}
k(\omega) = \frac{r^2+2r'^2-r r''}{(r^2+r'^2)^{3/2}} .
\end{equation}
The integral in (\ref{e33a}) for $\tilde{\sigma}_1$ is the length of 
the boundary $\vert\partial\Omega\vert$ for an arbitrary domain, i.e.
\begin{equation} \label{e36}
\tilde{\sigma}_1 = - \frac{2}{\sqrt{\pi}} 
\frac{\vert\partial\Omega\vert}{\vert\Omega\vert} .
\end{equation}
The substitution of the relations (\ref{e34}) and (\ref{e35}) into 
(\ref{e33b}) yields
\begin{equation} \label{e37}
\tilde{\sigma}_2 = \frac{\pi}{\vert\Omega\vert} +
\frac{1}{2\vert\Omega\vert} \int_0^{2\pi} {\rm d}\varphi
\frac{r'^2 - r r''}{r^2+r'^2} = \frac{\pi}{\vert\Omega\vert} ,
\end{equation}
since the integral 
\begin{equation} \label{e38} 
\int_0^{2\pi} {\rm d} \varphi \frac{r'^2 - r r''}{r^2+r'^2} = - \left[
\arctan \left( \frac{r'}{r} \right) \right]_0^{2\pi} = 0  
\end{equation}
for an arbitrary $r(\varphi)$.
Note the universal form of $\vert\Omega\vert \tilde{\sigma}_2 = \pi$ 
valid for any star-shaped domain.
The form of the higher coefficients depends on the particular shape
of the domain boundary.  

\subsection{Concave domain}
If the boundary possess some concave parts, these parts can be replaced 
by circles centered outside of the domain. 
Consequently, the function $\tau(s;\omega,\xi)$ in (\ref{e21}) is replaced
by the solution of the differential equation (\ref{e19}) for a thick annulus 
of the inner radius $R=R(\omega)$ and the outer radius going to infinity,
\begin{equation} \label{e39}
\tilde{\tau}(s,r)= \frac{1}{s^2} \left( 
1 - \frac{H_0^{(1)}(i s r)}{H_0^{(1)}(i s R)} \right) ,
\end{equation}
where $H_{\nu}^{(1)}(z)=J_{\nu}(z)+i Y_{\nu}(z)$ denotes the Hankel 
function\cite{GR}.
After the derivation with respect to $\xi=r-R$ at the boundary $\xi=0$, 
the ratio $i H_1^{(1)}(isR)/H_0^{(1)}(isR)$ can also be expressed using 
asymptotic $s\to\infty$ expansions\cite{GR}, with the result
\begin{equation} \label{e40}
\frac{\partial\tilde{\tau}(s;\xi)}{\partial\xi}\Big\vert_{\xi=0} =
\frac{1}{s} \left[ 1 + \frac{1}{2sR} - \frac{1}{8(sR)^2} + \frac{1}{8(sR)^3} 
\right. \left.
-\frac{25}{128(sR)^4} + \frac{13}{32(sR)^5} \cdots \right] . 
\end{equation}
In comparison with the formula (\ref{e32}) valid for convex domains, 
the signs are changed at the odd powers of $1/(sR)$. 
This means that the result (\ref{e32}) is applicable also to concave parts
of the domain boundary if we take the negative value of the outer curvature 
$-1/R(\omega)$.

\subsection{Comparison with Savo's results}
It is instructive to compare our coefficients $\{ \tilde{\sigma}_j \}$
of the large-$s$ expansion (\ref{e32}), expressed via the local curvature 
in equations (\ref{e33a})-(\ref{e33f}), with the exact ones 
derived previously by Savo\cite{Savo98a,Savo98b}.
Savo's results are expressed as integrals over the boundary $\partial\Omega$ 
of certain differential operators acting on the shortest distance 
$l({\bf r})$ of a point ${\bf r}$ from $\partial\Omega$.  
For two-dimensional domains, we were able to rewrite Savo's complicated 
formulas in terms of the curvature $k(\omega)$ and of their derivatives.
The exact results for the first four coefficients $\sigma_1,\sigma_2,\sigma_3$ 
and $\sigma_4$ coincide with our formulas (\ref{e33a})-(\ref{e33d}).  
In higher orders, we find
\begin{eqnarray} 
\sigma_5 & = & \frac{1}{240\sqrt{\pi}} \frac{1}{\vert\Omega\vert}
\int_{\partial\Omega} {\rm d}\omega \left\{ 25 k^4(\omega) - 
8[k'(\omega)]^2 \right\} , \label{e41a} \\
\sigma_6 & = & \frac{1}{192} \frac{1}{\vert\Omega\vert}
\int_{\partial\Omega} {\rm d}\omega
\left\{ 13 k^5(\omega) + 80 k(\omega) [k'(\omega)]^2
+ 47 k^2(\omega) k''(\omega) \right\} ; \label{e41b}
\end{eqnarray}
the expressions for the next coefficients are too complicated in Savo's 
format to handle with.
Comparing these exact results with the corresponding formulas (\ref{e33e}) 
and (\ref{e33f}) we see that the approximation by the local curvature
involves correctly the powers of the boundary curvature, but misses
derivatives of the curvature along the boundary.
We shall see in the next section that neglecting of curvature derivatives 
has only a small effect on the obtained result. 
The advantage of the approximation by the local curvature 
(\ref{e33a})-(\ref{e33f}) consists in its simplicity, numerical adequacy 
and availability for extremely high orders.

\section{Pad\'e interpolation}
This section deals with estimates of the low-energy part of the spectrum of 
the (minus) Dirichlet Laplacian.
The calculation is based on an interpolation between the large-$s$ and 
small-$s$ expansions of $\tau(s)$.
The convergence of the results is good especially for the lowest eigenvalue 
$\lambda_1$ which is chosen as a test of our calculations. 

For large $t$, the survival probability (\ref{e8}) can be approximated
by a sum of first few exponentials with the lowest eigenvalues, 
$S(t) \simeq \sum_{j=1}^{m} \gamma_j^2 e^{-\lambda_jt}$.
Hence, in the Laplace picture (\ref{e14}), we have for general domain
\begin{equation} \label{e42}
\tau(s) \simeq \sum_{j=1}^{m} \frac{\gamma_j^2}{s^2+\lambda_j}
= d_0 + d_1 s + d_2 s^2 + d_3 s^3 + d_4 s^4 + \cdots ,
\end{equation} 
where the even coefficients are given by
\begin{equation} \label{e43}
d_{2k} = (-1)^k \sum_{j=1}^{m} \frac{\gamma_j^2}{\lambda_j^k}
\end{equation}
and the odd coefficients vanish $d_{2k+1} = 0$.
The coefficients $\{ d_{2k}\}$ are not at one's disposal for general domains.
For simple domains like the disk or the ellipse, the even coefficients
$d_0,d_2,\ldots$ are, in principle, available to an arbitrary 
order\cite{Travenec10}.
Within the framework of the $m$-truncation (\ref{e42}), knowledge of
the first $2m$ nonzero coefficients $d_0,d_2,\ldots,d_{4m-2}$ is sufficient 
to get very good estimates of the $m$ lowest eigenvalues $\lambda_j$ and
the coefficients $\gamma_j^2$, by solving the corresponding set 
of $2m$ nonlinear equations.
For example, in the case of the unit disk with the exact value
$\lambda_1 = 5.783186\ldots$, we obtain $\lambda_1 = 5.784128\ldots$ for $m=1$,
$\lambda_1 = 5.783187\ldots$ for $m=2$, etc. 
The lowest eigenvalue can be found alternatively as
$\lim_{k\to\infty} d_{2k}/d_{2k+2} = -\lambda_1$.

In the previous section, we derived the large-$s$ expansion of $\tau(s)$ 
in the approximation by the local curvature:
\begin{equation} \label{e44}
\tilde{\tau}(s) = \frac{1}{s^2} + \sum_{j=1}^{\infty}
\frac{\Gamma(j/2+1)}{s^{j+2}} \tilde{\sigma}_j ,
\end{equation}
where the coefficients $\{ \tilde{\sigma}_j \}$ are given by 
(\ref{e33a})-(\ref{e33f}).
Our strategy is to find an interpolation scheme between this large-$s$
expansion and the small-$s$ expansion (\ref{e42}) of $\tau(s)$.
The lowest eigenvalue $\lambda_1$ will be deduced directly from
this interpolation scheme.

The interpolation method is based on the following steps:
\medskip

\noindent (i) {\bf Proposal of the Pad\'e interpolation:}
With respect to the first relation in equation (\ref{e42}), 
the small-$s$ behavior of $\tau(s)$ can be represented as a rational of 
polynomials $P_{2m-2}(s)/Q_{2m}(s)$ in even powers of $s$.
The large-$s$ expansion of $\tau(s)$ (\ref{e44}) contains also odd powers of 
$1/s$.
If we try to expand the first relation of (\ref{e42}) in $1/s$, the odd powers
would be missing.
To connect both expansions by one interpolation formula, we include
also odd powers of $s$ in the rational $P_{2m-2}(s)/Q_{2m}(s)$.
The interpolation function is then the Pad\'e approximant
\begin{equation} \label{e45}
\tau(s) = [n/n+2] \equiv \frac{p_0 + p_1 s + \cdots + p_{n-1}s^{n-1} + s^n}{
q_0 + q_1 s + \cdots + q_{n+1} s^{n+1} + s^{n+2}} ,
\end{equation}
where for simplification of the notation we use $n=2m-2$ 
(to obtain more data, also odd values of $n$ will be taken).
\medskip

\noindent (ii) {\bf Determination of Pad\'e parameters:}
To fit the $n$ unknown $p_j$ and $n+2$ unknown $q_j$ Pad\'e parameters, 
we expand the Pad\'e approximant (\ref{e45}) from both sides $s\to\infty$ 
and $s\to 0$, in the variables $1/s$ and $s$, respectively.
The expansion in $1/s$ has to reproduce the series (\ref{e44}) up to the term 
of the order $1/s^{n+4}$, implying in this way $n+2$ conditions for
the coefficients $p_j$ and $q_j$.
The expansion in $s$ has to respect the symmetry of the small-$s$
expansion (\ref{e42}), keeping only even powers of $s$ nonzero.
Therefore, we add $n$ conditions $d_1=0$, $d_3=0$, $\cdots$, $d_{2n-1}=0$.
We end up with the complete system of $2n+2$ nonlinear equations 
for $p_j$ and $q_j$. 
From all possible solutions we choose the real ones.
This scheme, which leads to a nonlinear set of equations, differs from 
the standard Pad\'e fitting method, which uses input series data from only 
one side and therefore leads to a linear set of equations.

\begin{table}
{Table 1. The coefficients of the small-$s$ expansion (46) and 
the lowest imaginary parts of poles of the Pad\'e approximant for the disk}
{\begin{tabular}{cccccc}
 \hline
 Pad\'e &$d_0$ &$d_2$ &$ d_4$ &$d_6$& Im[$s$]\\
 \hline
 \ [1/3]& 0.1743 & -0.07472 & 0.008383 & -0.00004709 & 1.756 \\
 \ [2/4]& 0.1475 & -0.03538 & 0.011763 & -0.00246499 & 1.940 \\
 \ [3/5]& 0.1378 & -0.02803 & 0.006355 & -0.00162755 & 2.074 \\
 \ [4/6]& 0.1331 & -0.02516 & 0.005090 & -0.00106198 & 2.178 \\
 \ [5/7]& 0.1306 & -0.02371 & 0.004541 & -0.00088446 & 2.252 \\
 \ [6/8]& 0.1290 & -0.02287 & 0.004246 & -0.00079754 & 2.299 \\
 \ [7/9]& 0.1289 & -0.02235 & 0.004067 & -0.00074735 & 2.328 \\
 \hline
 \ exact & 0.1250 & -0.02083 & 0.003581 & -0.00061849 & 2.405 \\
\end{tabular}}
\end{table}

\noindent (iii) {\bf Determination of $\lambda_1$ from Pad\'e approximant:} 
From the first relation of equation (\ref{e42}) we see that the poles of 
$\tau(s)$ occur in the complex plane at the points 
$\pm {\rm i}\sqrt{\lambda_j}$ $(j=1,2,\ldots)$.
The poles of the rational (\ref{e45}) are given by
\begin{equation} \label{e46}
q_0 + q_1 s + \cdots + q_{n+1} s^{n+1} + s^{n+2} = 0 .
\end{equation}
Since the coefficients $\{ q_j \}$ are real numbers, the roots
of this equation arise in complex conjugate pairs.
As $n\to\infty$, their (negative) real part should approach 0.
The lowest eigenvalue $\lambda_1$ is determined as the square of 
the pole closest to the origin.  
\medskip

As a test model, we consider the disk of radius $R=1$ for which 
the exact small-$s$ and large-$s$ expansions are given by the
formulas (\ref{e28}) and (\ref{e29}), respectively. 
We carried out numerical calculations up to $n=7$.
In Tab. 1, the first few nonzero coefficients of the small-$s$ expansion 
(\ref{e43}) are shown to converge very well to their exact values.
Fig. 2 depicts the poles of the Pad\'e approximant (\ref{e45}) which are
closest to the imaginary axis and converge to 
$s=\pm i\sqrt{\lambda_1}=\pm i z_{0,1}$ for $n\to\infty$.
It is seen that by increasing the Pad\'e order $n$, the (always negative) 
real part of the pole goes quickly to zero, as expected. 
The imaginary part approaches the exact value denoted by a cross, 
with the relative error 0.02 (after extrapolation).
The values of these imaginary parts of the poles are presented in Tab. 1.
In contrast to other numerical methods, like the harmonic inversion 
method\cite{Weibert05} or the finite elements method\cite{Son09},
our analytic treatment applies easily to general domains.

\vskip 0.6cm

\begin{figure}
\begin{center}
\includegraphics[scale=1]{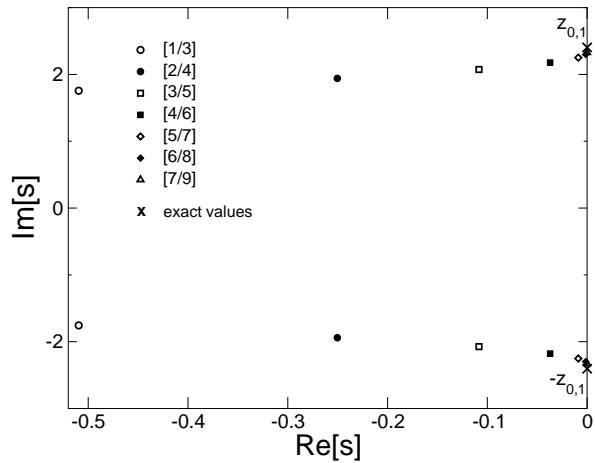}
\caption{The lowest poles of Pad\'e approximants $[n/n+2]$
$(n=1, ..., 7)$ for the disk.}
\end{center}
\end{figure}

\begin{figure}
\begin{center}
\includegraphics[scale=0.3]{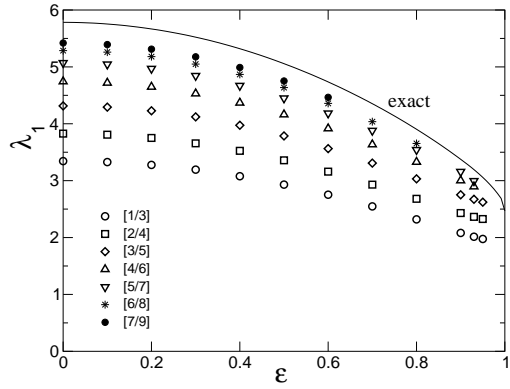}
\caption{Ellipse: the lowest eigenvalue $\lambda_1$ vs. eccentricity
$\varepsilon$ for Pad\'e approximants [n,n+2] $(n=1,2,\ldots,7)$.}
\end{center}
\end{figure}

For the elliptic domain with eccentricity $\varepsilon\in (0,1)$, 
we complete the integrations in (\ref{33}) in polar coordinates, in which
the boundary is given by $r(\varphi)=b/\sqrt{1-\epsilon^2\cos^2\varphi}$,
$b$ stands for the minor semiaxis.
The dependence of the function $\lambda_1(\varepsilon)$ on 
the Pad\'e order $n$ is presented in Fig. 3. 
We see that by increasing $n$ our estimates approach the exact 
plot\cite{Stegun}. 

The comparison of the results for the function $\lambda_1(\varepsilon)$, 
evaluated by using the large-$s$ coefficients in the approximation 
by the curvature (\ref{e33a})-(\ref{e33f}) and the exact 
ones\cite{Savo98a,Savo98b} [see (\ref{e41a}) and (\ref{e41b})],  
are presented in Fig. 4.
The results are almost identical up to very large eccentricity 
$\varepsilon\sim 0.9$, the difference remains small also beyond 
this point.

\vskip 0.6cm

\begin{figure}
\begin{center}
\includegraphics[scale=0.3]{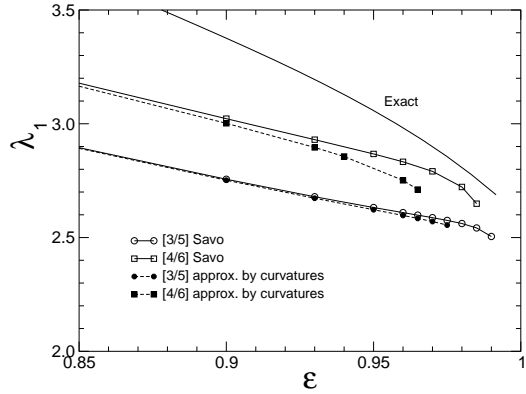}
\caption{Ellipse: $\lambda_1$ vs. eccentricity $\varepsilon$ for 
Pad\'e orders $n=3,4$; the comparison of estimates obtained by using 
the approximation by the curvature (full symbols) and Savo's results 
(\ref{e41a}), (\ref{e41b}) (empty symbols).}
\end{center}
\end{figure}

\section{Conclusion}
In this work, we studied the survival probability $S(t)$ of a particle
diffusing in a finite domain with a smooth absorbing boundary. 
Its short time expansion in $\sqrt{t}$ (\ref{e9}) is determined by 
the boundary characteristics, whilst its long time exponential decay 
is governed by the lowest eigenvalue of the minus Laplacian defined 
in this domain. 
Thus this quantity enables us to study the impact of the boundary shape 
on the lowest eigenvalue $\lambda_1$.

First we proposed a simple method which allows us to calculate 
the coefficients $\sigma_j$ of the short time expansion of $S(t)$ 
up to high orders.
We approximated the boundary locally by circles of the radius inverse to
the local curvature of the boundary. 
Using the known exact solution for the disk domain, we expressed 
the coefficients $\sigma_j$ in terms of powers of the curvature 
integrated along the circumference of the domain, see equations
(\ref{e33a})-(\ref{e33f}). 
This treatment gives exact results up to $\sigma_4$.
In higher orders, the exact coefficients of Savo\cite{Savo98a,Savo98b} 
contain also derivatives of the curvature, which are neglected in 
the present approximation.
While Savo's results can be obtained in practice only in the first
few orders, our formulas are at disposal for virtually any order.

Next, we showed that there exists an interpolation scheme connecting the
short time expansion of $S(t)$ with the lowest eigenvalues of 
the Dirichlet Laplacian.
It is represented by the $[n/n+2]$ Pad\'e approximation (\ref{e45}) of 
the Laplace transform $\tau(s)$ of the survival probability $S(t)$ (\ref{e14}). 
If the coefficients $p_i$, $q_j$ are fixed in a specific way, the poles 
closest to the imaginary axis tend to $\pm {\rm i}\sqrt{\lambda_1}$.
Extrapolation of our data for $n\le 7$, calculated for the disk and the 
elliptic domains, exhibits precision within 2\% of the exact value.

Of course, there may exist other interpolation schemes which provide
faster convergence to the exact value of $\lambda_1$.
Our present aim was mainly to point out the existence of such an algorithm.
Another task is to extend the calculation of $S(t)$ to domains with 
cusps (with piecewise smooth boundaries), or to higher-dimensional domains. 
These are open problems for future.

\vskip 0.6truecm

\section*{Acknowledgments}
This work was supported by the Grants VEGA No. 2/0113/2010 and CE-SAS QUTE.

\end{document}